%%%%%%%%%%%%%%%%%%%%%%%%%%%%%%%%%%%%%%%%%%%%%%%%%%%%%%%%%%%%%%%%%%%%%%%%%%%%%
%%%%%%%%%%%%%%%%%%%%%%%%%  DEFINIZIONI DI FONT    %%%%%%%%%%%%%%%%%%%%%%%%%%%
%%%%%%%%%%%%%%%%%%%%%%%%%%%%%%%%%%%%%%%%%%%%%%%%%%%%%%%%%%%%%%%%%%%%%%%%%%%%%
\font\tenmib=cmmib10
\font\sevenmib=cmmib10 scaled 800
\font\titolo=cmbx12
\font\titolone=cmbx10 scaled\magstep 2

\font\cs=cmcsc10
\font\sc=cmcsc10
\font\css=cmcsc8

\font\ninerm=cmr9
\font\ottorm=cmr8
\textfont5=\tenmib\scriptfont5=\sevenmib\scriptscriptfont5=\fivei

\font\msytw=msbm9 scaled\magstep1

\font\msytwww=msbm7 scaled\magstep1
\font\msytwwww=msbm5 scaled\magstep1
\font\msytwwwww=msbm4 scaled\magstep1
\font\msytwwwwww=msbm3 scaled\magstep1
\font\indbf=cmbx10 scaled\magstep2

\font\ottorm=cmr8\font\ottoi=cmmi8\font\ottosy=cmsy8%
\font\ottobf=cmbx8\font\ottott=cmtt8%
\font\ottocss=cmcsc8%
\font\ottosl=cmsl8\font\ottoit=cmti8%
\font\sixrm=cmr6\font\sixbf=cmbx6\font\sixi=cmmi6\font\sixsy=cmsy6%
\font\fiverm=cmr5\font\fivesy=cmsy5\font\fivei=cmmi5\font\fivebf=cmbx5%

\def\ottopunti{\def\rm{\fam0\ottorm}%
\textfont0=\ottorm\scriptfont0=\sixrm\scriptscriptfont0=\fiverm%
\textfont1=\ottoi\scriptfont1=\sixi\scriptscriptfont1=\fivei%
\textfont2=\ottosy\scriptfont2=\sixsy\scriptscriptfont2=\fivesy%
\textfont3=\tenex\scriptfont3=\tenex\scriptscriptfont3=\tenex%
\textfont4=\ottocss\scriptfont4=\sc\scriptscriptfont4=\sc%
%\scriptfont4=\ottocss\scriptscriptfont4=\ottocss%
%\textfont5=\tenmib\scriptfont5=\sevenmib\scriptscriptfont5=\fivei%
\textfont\itfam=\ottoit\def\it{\fam\itfam\ottoit}%
\textfont\slfam=\ottosl\def\sl{\fam\slfam\ottosl}%
\textfont\ttfam=\ottott\def\tt{\fam\ttfam\ottott}%
\textfont\bffam=\ottobf\scriptfont\bffam=\sixbf%
\scriptscriptfont\bffam=\fivebf\def\bf{\fam\bffam\ottobf}%
%\tt\ttglue=.5em plus.25em minus.15em%
\setbox\strutbox=\hbox{\vrule height7pt depth2pt width0pt}%
\normalbaselineskip=9pt\let\sc=\sixrm\normalbaselines\rm}
%

%%%%%%%%%%%%%%%%%%%%%%%%%%%%%%%%%%%%%%%%%%%%%%%%%%%%%%%%%%%%%%%%%%%%%%%%%%%%%
%%%%%%%%%   RIFERIMENTI SIMBOLICI A FORMULE, PARAGRAFI E FIGURE    %%%%%%%%%%
%%%%%%%%%%%%%%%%%%%%%%%%%%%%%%%%%%%%%%%%%%%%%%%%%%%%%%%%%%%%%%%%%%%%%%%%%%%%%
%
% Ogni paragrafo deve iniziare con il comando \section(#1,#2), dove #1
% e' il simbolo associato al paragrafo e #2 e' il titolo. Per le
% appendici bisogna pero' usare \appendix(#1,#2).
%
% Se nel titolo compaiono riferimenti ad altri simboli, questi vanno
% racchiusi fra parentesi graffe, per es. {\equ(1.2)}; in caso contrario
% si provoca un errore.
%
% Ogni sottoparagrafo deve iniziare con il comando \sub(#1) o \asub(#1),
% nelle appendici.
%
% I riferimenti a paragrafi e sottoparagrafi si realizzano con il comando
% \sec(#1), che produce il numero effettivo preceduto dal simbolo di
% paragrafo, o \secc(#1), che produce solo il numero (serve nel caso si
% faccia riferimento ad un sottoparagrafo, che e' un Lemma, un Teorema o
% altro oggetto suscettibile di una denominazione speciale).
%
% Le formule sono contrassegnate con \Eq(#1), eccetto che all'interno
% del comando \eqalignno, dove si deve usare \eq(#1). Nelle appendici
% i comandi corrispondenti sono \Eqa(#1) e \eqa(#1).
% I riferimenti alle formule si realizzano con \equ(#1).
%
% La numerazione delle figure utilizza il comando \eqg(#1), per
% contrassegnarle, e \graf(#1) per citarle.
%

\global\newcount\numsec\global\newcount\numapp
\global\newcount\numfor\global\newcount\numfig
\global\newcount\numsub
\numsec=0\numapp=0\numfig=0
\def\veroparagrafo{\number\numsec}\def\veraformula{\number\numfor}
\def\veraappendice{\number\numapp}\def\verasub{\number\numsub}
\def\verafigura{\number\numfig}

\def\section(#1,#2){\advance\numsec by 1\numfor=1\numsub=1\numfig=1%
\SIA p,#1,{\veroparagrafo} %
\write15{\string\Fp (#1){\secc(#1)}}%
\write16{ sec. #1 ==> \secc(#1)  }%
\hbox to \hsize{\titolo\hfill \number\numsec. #2\hfill%
\expandafter{\alato(sec. #1)}}\*}

\def\appendix(#1,#2){\advance\numapp by 1\numfor=1\numsub=1\numfig=1%
\SIA p,#1,{A\veraappendice} %
\write15{\string\Fp (#1){\secc(#1)}}%
\write16{ app. #1 ==> \secc(#1)  }%
\hbox to \hsize{\titolo\hfill Appendix A\number\numapp. #2\hfill%
\expandafter{\alato(app. #1)}}\*}

\def\senondefinito#1{\expandafter\ifx\csname#1\endcsname\relax}

\def\SIA #1,#2,#3 {\senondefinito{#1#2}%
\expandafter\xdef\csname #1#2\endcsname{#3}\else
\write16{???? ma #1#2 e' gia' stato definito !!!!} \fi}

\def \Fe(#1)#2{\SIA fe,#1,#2 }
\def \Fp(#1)#2{\SIA fp,#1,#2 }
\def \Fg(#1)#2{\SIA fg,#1,#2 }

\def\etichetta(#1){(\veroparagrafo.\veraformula)%
\SIA e,#1,(\veroparagrafo.\veraformula) %
\global\advance\numfor by 1%
\write15{\string\Fe (#1){\equ(#1)}}%
\write16{ EQ #1 ==> \equ(#1)  }}

\def\etichettaa(#1){(A\veraappendice.\veraformula)%
\SIA e,#1,(A\veraappendice.\veraformula) %
\global\advance\numfor by 1%
\write15{\string\Fe (#1){\equ(#1)}}%
\write16{ EQ #1 ==> \equ(#1) }}

\def\getichetta(#1){\veroparagrafo.\verafigura%
\SIA g,#1,{\veroparagrafo.\verafigura} %
\global\advance\numfig by 1%
\write15{\string\Fg (#1){\graf(#1)}}%
\write16{ Fig. #1 ==> \graf(#1) }}

\def\getichettaa(#1){A\veraappendice.\verafigura%
\SIA g,#1,{A\veraappendice.\verafigura} %
\global\advance\numfig by 1%
\write15{\string\Fg (#1){\graf(#1)}}%
\write16{ Fig. #1 ==> \graf(#1) }}

\def\etichettap(#1){\veroparagrafo.\verasub%
\SIA p,#1,{\veroparagrafo.\verasub} %
\global\advance\numsub by 1%
\write15{\string\Fp (#1){\secc(#1)}}%
\write16{ par #1 ==> \secc(#1)  }}

\def\etichettapa(#1){A\veraappendice.\verasub%
\SIA p,#1,{A\veraappendice.\verasub} %
\global\advance\numsub by 1%
\write15{\string\Fp (#1){\secc(#1)}}%
\write16{ par #1 ==> \secc(#1)  }}

\def\Eq(#1){\eqno{\etichetta(#1)\alato(#1)}}
\def\eq(#1){\etichetta(#1)\alato(#1)}
\def\Eqa(#1){\eqno{\etichettaa(#1)\alato(#1)}}
\def\eqa(#1){\etichettaa(#1)\alato(#1)}
\def\eqg(#1){\getichetta(#1)\alato(fig. #1)}
\def\sub(#1){\0\palato(p. #1){\bf \etichettap(#1).}}
\def\asub(#1){\0\palato(p. #1){\bf \etichettapa(#1).}}

\def\equv(#1){\senondefinito{fe#1}$\clubsuit$#1%
\write16{eq. #1 non e' (ancora) definita}%
\else\csname fe#1\endcsname\fi}
\def\grafv(#1){\senondefinito{fg#1}$\clubsuit$#1%
\write16{fig. #1 non e' (ancora) definito}%
\else\csname fg#1\endcsname\fi}
\def\secv(#1){\senondefinito{fp#1}$\clubsuit$#1%
\write16{par. #1 non e' (ancora) definito}%
\else\csname fp#1\endcsname\fi}

\def\equ(#1){\senondefinito{e#1}\equv(#1)\else\csname e#1\endcsname\fi}
\def\graf(#1){\senondefinito{g#1}\grafv(#1)\else\csname g#1\endcsname\fi}
\def\figura(#1){{\css Figure} \getichetta(#1)}
\def\figuraa(#1){{\css Figure} \getichettaa(#1)}
\def\secc(#1){\senondefinito{p#1}\secv(#1)\else\csname p#1\endcsname\fi}
\def\sec(#1){{\S\secc(#1)}}
\def\refe(#1){{[\secc(#1)]}}

\def\BOZZA{\bz=1
\def\alato(##1){\rlap{\kern-\hsize\kern-1.2truecm{$\scriptstyle##1$}}}
\def\palato(##1){\rlap{\kern-1.2truecm{$\scriptstyle##1$}}}
}

\def\alato(#1){}
\def\galato(#1){}
\def\palato(#1){}

%%%%%%%%%%%%%%%%%%%%%%%%%%%%%%%%%%%%%%%%%%%%%%%%%%%%%%%%%%%%%%%%%%%%%%%%%%%%%
%%%%%%%%%%%%%%%%%%%%      DATA E PIE' DI PAGINA        %%%%%%%%%%%%%%%%%%%%%%
%%%%%%%%%%%%%%%%%%%%%%%%%%%%%%%%%%%%%%%%%%%%%%%%%%%%%%%%%%%%%%%%%%%%%%%%%%%%%

{\count255=\time\divide\count255 by 60 \xdef\hourmin{\number\count255}
        \multiply\count255 by-60\advance\count255 by\time
   \xdef\hourmin{\hourmin:\ifnum\count255<10 0\fi\the\count255}}

\def\oramin{\hourmin }

\def\data{\number\day/\ifcase\month\or gennaio \or febbraio \or marzo \or
aprile \or maggio \or giugno \or luglio \or agosto \or settembre
\or ottobre \or novembre \or dicembre \fi/\number\year;\ \oramin}
%\setbox200\hbox{$\scriptscriptstyle \data $}
\footline={\rlap{\hbox{\copy200}}\tenrm\hss \number\pageno\hss}

%%%%%%%%%%%%%%%%%%%%%%%%%%%%%%%%%%%%%%%%%%%%%%%%%%%%%%%%%%%%%%%%%%%%%%%%%%%%%
%%%%%%%%%%%%%%%      INSERIMENTO FIGURE ( se si usa DVIPS )    %%%%%%%%%%%%%%
%%%%%%%%%%%%%%%%%%%%%%%%%%%%%%%%%%%%%%%%%%%%%%%%%%%%%%%%%%%%%%%%%%%%%%%%%%%%%

\newcount\driver %\driver=1 %dvips
\newdimen\xshift \newdimen\xwidth
\def\ins#1#2#3{\vbox to0pt{\kern-#2 \hbox{\kern#1
#3}\vss}\nointerlineskip}

\def\insertplot#1#2#3#4#5{\par%
\xwidth=#1 \xshift=\hsize \advance\xshift by-\xwidth \divide\xshift by 2%
\yshift=#2 \divide\yshift by 2%
\line{\hskip\xshift \vbox to #2{\vfil%
\ifnum\driver=0 #3
%
\special{ps: plotfile #4.eps}
\fi
\ifnum\driver=1 #3 \includegraphics{#4.eps}\fi
\ifnum\driver=2 #3
\ifnum\mgnf=0\special{#4.eps 1. 1. scale} \fi
\ifnum\mgnf=1\special{#4.eps 1.2 1.2 scale}\fi
\fi }\hfill \raise\yshift\hbox{#5}}}

\def\insertplotttt#1#2#3#4{\par%
\xwidth=#1 \xshift=\hsize \advance\xshift by-\xwidth \divide\xshift by 2%
\yshift=#2 \divide\yshift by 2%
\line{\hskip\xshift \vbox to #2{\vfil%
\ifnum\driver=0 #3
%
\special{ps: plotfile #4.eps}
\fi
\ifnum\driver=1 #3 \includegraphics{#4.eps}\fi
\ifnum\driver=2 #3
\ifnum\mgnf=0\special{#4.eps 1. 1. scale} \fi
\ifnum\mgnf=1\special{#4.eps 1.2 1.2 scale}\fi
\fi }\hfill}}

\def\insertplott#1#2#3{\par%
\xwidth=#1 \xshift=\hsize \advance\xshift by-\xwidth \divide\xshift by 2%
\yshift=#2 \divide\yshift by 2%
\line{\hskip\xshift \vbox to #2{\vfil%
\ifnum\driver=0
%
\special{ps: plotfile #3.eps}
 \fi
\ifnum\driver=1 \includegraphics{#3.eps}\fi
\ifnum\driver=2
\ifnum\mgnf=0\special{#3.eps 1. 1. scale} \fi 
\ifnum\mgnf=1\special{#3.eps 1.2 1.2 scale}\fi
\fi }\hfill}}

\newdimen\xshift \newdimen\xwidth \newdimen\yshift
\def\eqfig#1#2#3#4#5{
\par\xwidth=#1 \xshift=\hsize \advance\xshift
by-\xwidth \divide\xshift by 2
\yshift=#2 \divide\yshift by 2
\line{\hglue\xshift \vbox to #2{\vfil
\ifnum\driver=0 #3
%
\special{ps: plotfile #4.ps}
\fi
\ifnum\driver=1 #3 \includegraphics{#4.ps}\fi
\ifnum\driver=2 #3 \special{
\ifnum\mgnf=0 #4.ps 1. 1. scale \fi
\ifnum\mgnf=1 #4.ps 1.2 1.2 scale\fi}
\fi}\hfill\raise\yshift\hbox{#5}}}

\let\a=\alpha \let\b=\beta  \let\g=\gamma  \let\d=\delta \let\e=\varepsilon
  \let\h=\eta   \let\th=\theta \let\k=\kappa 
\let\m=\mu    \let\n=\nu         \let\p=\pi    
\let\s=\sigma \let\t=\tau   \let\f=\varphi 
   \let\o=\omega
   \let\Th=\Theta 
     \let\F=\Phi

\def\\{\hfill\break} \let\==\equiv

\let\io=\infty 

\let\0=\noindent

\let\dpr=\partial

\def\tende#1{\,\vtop{\ialign{##\crcr\rightarrowfill\crcr
 \noalign{\kern-1pt\nointerlineskip}
 \hskip3.pt${\scriptstyle #1}$\hskip3.pt\crcr}}\,}
\def\otto{\,{\kern-1.truept\leftarrow\kern-5.truept\to\kern-1.truept}\,}

\def\EE{{\cal E}}\def\MM{{\cal M}} \def\VV{{\cal V}}

\def\LL{{\cal L}} 
\def\DD{{\cal D}}\def\AAA{{\cal A}}

\def\T#1{{#1_{\kern-3pt\lower7pt\hbox{$\widetilde{}$}}\kern3pt}}
\def\VVV#1{{\underline #1}_{\kern-3pt
\lower7pt\hbox{$\widetilde{}$}}\kern3pt\,}
\def\W#1{#1_{\kern-3pt\lower7.5pt\hbox{$\widetilde{}$}}\kern2pt\,}

\def\lis{\overline}

\def\indica{\leaders \hbox to 0.5cm{\hss.\hss}\hfill}
\def\guida{\leaders\hbox to 1em{\hss.\hss}\hfill}

   \def\V0{{\bf 0}}

\def\oo{{\omega}}
%%%%%%%%%%%%%%%%%%%%%%%%%%%%%%%%%%%%%%%%%%%%%%%%%%%%%%%%%%%%%%%%%%%%%%%%%%%%%
%%%%%%%%%%%%%%%%%    LETTERE GRECHE E LATINE IN NERETTO     %%%%%%%%%%%%%%%%%
%%%%%%%%%%%%%%%%%%%%%%%%%%%%%%%%%%%%%%%%%%%%%%%%%%%%%%%%%%%%%%%%%%%%%%%%%%%%%

% lettere greche e latine in neretto italico - pag.430 del manuale
\mathchardef\aa   = "050B
\mathchardef\bb   = "050C
\mathchardef\xxx  = "0518
\mathchardef\hhh  = "0511
\mathchardef\zzzzz= "0510
\mathchardef\oo   = "0521
\mathchardef\tt   = "051C
\mathchardef\lll  = "0515
\mathchardef\mmm  = "0516
\mathchardef\Dp   = "0540
\mathchardef\H    = "0548
\mathchardef\FFF  = "0546
\mathchardef\ppp  = "0570
\mathchardef\nnn  = "0517
\mathchardef\pps  = "0520
\mathchardef\XXX  = "0504
\mathchardef\FFF  = "0508
\mathchardef\nnnnn= "056E
\mathchardef\minore  = "053C
\mathchardef\maggiore  = "053E

\def\to{\rightarrow}

\let\ciao=\bye
\def\qed{\hfill\raise1pt\hbox{\vrule height5pt width5pt depth0pt}}

 \def\Val{{\rm Val}}
\def\indic{\hbox{\raise-2pt \hbox{\indbf 1}}}

\def\RRR{\hbox{\msytw R}}

\def\NNN{\hbox{\msytw N}} 
\def\nnn{\hbox{\msytwww N}}
\def\ZZZ{\hbox{\msytw Z}} 
\def\zzz{\hbox{\msytwww Z}}

\def\vvv{\hbox{\msytwwww V}} \def\vvvv{\hbox{\msytwwwww V}}
\def\www{\hbox{\msytwwww W}} \def\wwww{\hbox{\msytwwwww W}}
\def\wwwww{\hbox{\msytwwwwww W}}

%%%%%%%%%%%%%%%%%%%%%%%%%%%%%%%%%%%%%%%%%%%%%%%%%%%%%%%%%%%%%%%%%%%%%%%%%%%%%
%%%%%%%%%%%%%%%%%    FORMATO PAGINA E LETTURA FILE .AUX     %%%%%%%%%%%%%%%%%
%%%%%%%%%%%%%%%%%%%%%%%%%%%%%%%%%%%%%%%%%%%%%%%%%%%%%%%%%%%%%%%%%%%%%%%%%%%%%

\newcount\mgnf  %ingrandimento
\mgnf=0 

\ifnum\mgnf=0
\def\openone{\leavevmode\hbox{\ninerm 1\kern-3.3pt\tenrm1}}%
\def\*{\vglue0.3truecm}\fi
\ifnum\mgnf=1
\def\openone{\leavevmode\hbox{\ninerm 1\kern-3.63pt\tenrm1}}%
\def\*{\vglue0.5truecm}\fi

%%%%%%%%%%%%%%%%%%%%%%%%%%%%%%%%%%%%%%%%%%%%%%%%%%%%%%%%%%%%%%%%%%%%%%%%%%%%%
%%%%%%%%%%%%%%%%%%%%%%         BIBBLIOGRAFIA        %%%%%%%%%%%%%%%%%%%%%%%%%
%%%%%%%%%%%%%%%%%%%%%%%%%%%%%%%%%%%%%%%%%%%%%%%%%%%%%%%%%%%%%%%%%%%%%%%%%%%%%

\newcount\tipobib\newcount\bz\bz=0\newcount\aux\aux=1
\newdimen\bibskip\newdimen\maxit\maxit=0pt

%\openin14=\jobname.aux \ifeof14 \relax \else
%\input \jobname.aux \closein14 \fi
%\openout15=\jobname.aux

\tipobib=0
\def\9#1{\ifnum\aux=1#1\else\relax\fi}

\newwrite\bib
\immediate\openout\bib=\jobname.bib
\global\newcount\bibex
\bibex=0
\def\verabib{\number\bibex}

\ifnum\tipobib=0
\def\cita#1{\expandafter\ifx\csname c#1\endcsname\relax
\hbox{$\clubsuit$}#1\write16{Manca #1 !}%
\else\csname c#1\endcsname\fi}
\def\rife#1#2#3{\immediate\write\bib{\string\raf{#2}{#3}{#1}}
\immediate\write15{\string\C(#1){[#2]}}
\setbox199=\hbox{#2}\ifnum\maxit < \wd199 \maxit=\wd199\fi}
\else
\def\cita#1{%
\expandafter\ifx\csname d#1\endcsname\relax%
\expandafter\ifx\csname c#1\endcsname\relax%
\hbox{$\clubsuit$}#1\write16{Manca #1 !}%
\else\probib(ref. numero )(#1)%
\csname c#1\endcsname%
\fi\else\csname d#1\endcsname\fi}%
\def\rife#1#2#3{\immediate\write15{\string\Cp(#1){%
\string\immediate\string\write\string\bib{\string\string\string\raf%
{\string\verabib}{#3}{#1}}%
\string\Cn(#1){[\string\verabib]}%
\string\CCc(#1)%
}}}%
\fi

\def\Cn(#1)#2{\expandafter\xdef\csname d#1\endcsname{#2}}
\def\CCc(#1){\csname d#1\endcsname}
\def\probib(#1)(#2){\global\advance\bibex+1%
\9{\immediate\write16{#1\verabib => #2}}%
}

\def\C(#1)#2{\SIA c,#1,{#2}}
\def\Cp(#1)#2{\SIAnx c,#1,{#2}}

\def\SIAnx #1,#2,#3 {\senondefinito{#1#2}%
\expandafter\def\csname#1#2\endcsname{#3}\else%
\write16{???? ma #1,#2 e' gia' stato definito !!!!}\fi}

\bibskip=10truept
\def\hboxto{\hbox to}

\catcode`\{=12\catcode`\}=12
\catcode`\<=1\catcode`\>=2
\immediate\write\bib<
        \string\halign{\string\hboxto \string\maxit% 
        {\string #\string\hfill}&% 
        \string\vtop{\string\parindent=0pt\string\advance\string\hsize% 
        by -1.55truecm%
        \string#\string\vskip \bibskip
        }\string\cr%
>
\catcode`\{=1\catcode`\}=2
\catcode`\<=12\catcode`\>=12

\def\raf#1#2#3{\ifnum \bz=0 [#1]&#2 \cr\else
\llap{${}_{\rm #3}$}[#1]&#2\cr\fi}

\newread\bibin

\catcode`\{=12\catcode`\}=12
\catcode`\<=1\catcode`\>=2
\def\chiudibib<
\catcode`\{=12\catcode`\}=12
\catcode`\<=1\catcode`\>=2
\immediate\write\bib<}>
\catcode`\{=1\catcode`\}=2
\catcode`\<=12\catcode`\>=12
>
\catcode`\{=1\catcode`\}=2
\catcode`\<=12\catcode`\>=12

\def\makebiblio{
\ifnum\tipobib=0
\advance \maxit by 10pt
\else
\maxit=1.truecm
\fi
\chiudibib
\immediate \closeout\bib
\openin\bibin=\jobname.bib
\ifeof\bibin\relax\else
\raggedbottom
\input \jobname.bib
\fi
}

\openin13=#1.aux \ifeof13 \relax \else
\input #1.aux \closein13\fi
\openin14=\jobname.aux \ifeof14 \relax \else
\input \jobname.aux \closein14 \fi
\immediate\openout15=\jobname.aux

\def\biblio{\*\*\centerline{\titolo References}\*\nobreak\makebiblio}

%%%%%%%%%%%%%%%%%%%%%%%%%%%%%%%%%%%%%%%%%%%%%%%%%%%%%%%%%%%%%%%%%%%%%%%%%%%%%
%%%%%%%%%%%%%%%%%%%%%%%%         CORNICE              %%%%%%%%%%%%%%%%%%%%%%%
%%%%%%%%%%%%%%%%%%%%%%%%%%%%%%%%%%%%%%%%%%%%%%%%%%%%%%%%%%%%%%%%%%%%%%%%%%%%%

\ifnum\mgnf=0
   \magnification=\magstep0
   \hsize=14.truecm\vsize=22.0truecm\voffset0.truecm\hoffset.2truecm
   \parindent=0.3cm\baselineskip=0.43cm\fi
\ifnum\mgnf=1
   \magnification=\magstep1\hoffset=0.truecm
   \hsize=14.truecm\vsize=22.0truecm
   \baselineskip=18truept plus0.1pt minus0.1pt \parindent=0.9truecm
   \lineskip=0.5truecm\lineskiplimit=0.1pt      \parskip=0.1pt plus1pt\fi

%%% IMPOSTAZIONI PARAMETRI E LETTURA FILE .AUX

\mgnf=0   %ingrandimento
\driver=1 %figure
\openin14=\jobname.aux \ifeof14 \relax \else
\input \jobname.aux \closein14 \fi
\openout15=\jobname.aux

%%%%%%%%%%%%%%%%%%%%%%%%%%%%%%%%%%%%%%%%%%%%%%%%%%%%%%%%%%%%%%%%%%%%%%%%%%
%%%%%%%%%%%%%%%%%%%%%%%%%%%%%%%%%%%%%%%%%%%%%%%%%%%%%%%%%%%%%%%%%%%%%%%%%%
%%%%%%%%%%%%%%%%%%%%%%%%%%%%%%%%%%%%%%%%%%%%%%%%%%%%%%%%%%%%%%%%%%%%%%%%%%
%                                                                        %
%                              ARTICOLO                                  %
%                                                                        %
%%%%%%%%%%%%%%%%%%%%%%%%%%%%%%%%%%%%%%%%%%%%%%%%%%%%%%%%%%%%%%%%%%%%%%%%%%
%%%%%%%%%%%%%%%%%%%%%%%%%%%%%%%%%%%%%%%%%%%%%%%%%%%%%%%%%%%%%%%%%%%%%%%%%%
%%%%%%%%%%%%%%%%%%%%%%%%%%%%%%%%%%%%%%%%%%%%%%%%%%%%%%%%%%%%%%%%%%%%%%%%%%

%\BOZZA
%\headline{{\hss {\rm DRAFT - Not for circulation}}}

%\null\vskip1truecm

%%%%%%%%%%%%%%%%%%%%%%%%%%%%%%%%%%%%%%%%%%%%%%%%%%%%%%%%%%%%%%%%%%%%%%%%%%
%%%%%%%%%%%%%%%%%               TITOLO                  %%%%%%%%%%%%%%%%%%
%%%%%%%%%%%%%%%%%%%%%%%%%%%%%%%%%%%%%%%%%%%%%%%%%%%%%%%%%%%%%%%%%%%%%%%%%%
%%%%%%%%%%%%%%%%%%%%%%%%%%%%%%%%%%%%%%%%%%%%%%%%%%%%%%%%%%%%%%%%%%%%%%%%%%
\centerline{\titolone Conservation of resonant periodic solutions for the}
\vskip.1cm
\centerline{\titolone one-dimensional nonlinear Schr\"odinger equation}
%%%%%%%%%%%%%%%%%%%%%%%%%%%%%%%%%%%%%%%%%%%%%%%%%%%%%%%%%%%%%%%%%%%%%%%%%%
%%%%%%%%%%%%%%%%%%%%%%%%%%%%%%%%%%%%%%%%%%%%%%%%%%%%%%%%%%%%%%%%%%%%%%%%%%

%%%%%%%%%%%%%%%%%%%%%%%%%%%%%%%%%%%%%%%%%%%%%%%%%%%%%%%%%%%%%%%%%%%%%%%%%%
%%%%%%%%%%%%%%%%%               AUTORI                  %%%%%%%%%%%%%%%%%%
%%%%%%%%%%%%%%%%%%%%%%%%%%%%%%%%%%%%%%%%%%%%%%%%%%%%%%%%%%%%%%%%%%%%%%%%%%
\vskip1.truecm \centerline{{\titolo
Guido Gentile$^{\dagger}$ and Michela Procesi$^{\star}$}}
\vskip.2truecm
\centerline{{}$^{\dagger}$ Dipartimento di Matematica,
Universit\`a di Roma Tre, Roma, I-00146}
\centerline{{}$^{\star}$ SISSA, Trieste, I-34014}
\vskip1.truecm
%%%%%%%%%%%%%%%%%%%%%%%%%%%%%%%%%%%%%%%%%%%%%%%%%%%%%%%%%%%%%%%%%%%%%%%%%%

%%%%%%%%%%%%%%%%%%%%%%%%%%%%%%%%%%%%%%%%%%%%%%%%%%%%%%%%%%%%%%%%%%%%%%%%%%
%%%%%%%%%%%%%%%%%             RIASSUNTO                 %%%%%%%%%%%%%%%%%%
%%%%%%%%%%%%%%%%%%%%%%%%%%%%%%%%%%%%%%%%%%%%%%%%%%%%%%%%%%%%%%%%%%%%%%%%%%
\line{\vtop{
\line{\hskip1.1truecm\vbox{\advance \hsize by -2.3 truecm
\0{\cs Abstract.}
{\it 
We consider the one-dimensional nonlinear Schr\"odinger equation
with Dirichlet boundary conditions in the fully
resonant case (absence of the zero-mass term).
We investigate conservation of small amplitude periodic-solutions
for a large set measure of frequencies. In particular we show that
there are infinitely many periodic solutions which continue
the linear ones involving an arbitrary number of resonant modes,
provided the corresponding frequencies are large enough
and close enough to each other (wave packets with large wave number).
}} \hfill} }}
%%%%%%%%%%%%%%%%%%%%%%%%%%%%%%%%%%%%%%%%%%%%%%%%%%%%%%%%%%%%%%%%%%%%%%%%%%
\vskip0.5truecm

%%%%%%%%%%%%%%%%%%%%%%%%%%%%%%%%%%%%%%%%%%%%%%%%%%%%%%%%%%%%%%%%%%%%%%%%%%
%%%%%%%%%%%%%%%%%%%%%%%%%%%%%%%%%%%%%%%%%%%%%%%%%%%%%%%%%%%%%%%%%%%%%%%%%%
\*\*
%%%%%%%%%%%%%%%%%%%%%%%%%%%%%%%%%%%%%%%%%%%%%%%%%%%%%%%%%%%%%%%%%%%%%%%%%%
%%%%%%%%%%%%%%%%%%%%%%%%%%%%%%%%%%%%%%%%%%%%%%%%%%%%%%%%%%%%%%%%%%%%%%%%%%
\section(1,Introduction and set-up)
%%%%%%%%%%%%%%%%%%%%%%%%%%%%%%%%%%%%%%%%%%%%%%%%%%%%%%%%%%%%%%%%%%%%%%%%%%
%%%%%%%%%%%%%%%%%%%%%%%%%%%%%%%%%%%%%%%%%%%%%%%%%%%%%%%%%%%%%%%%%%%%%%%%%%

\0We consider the {\it nonlinear Sch\"rodinger equation}
in $d=1$ given by
$$ \cases{ 
- i u_{t} + u_{xx} = \f(|u|^{2})u , & \cr
u(0,t) = u(\p,t) = 0 , & \cr}
\Eq(1.1)$$
where $\f(x)$ is any analytic function
$\f(x) = \F x + O(x^{2})$ with $\F\neq0$.
We shall consider the problem of existence of periodic solutions for
\equ(1.1), and we shall show how suitably adapting
the techniques in Ref.~\cita{GMP} we can solve the problem.

Existence of periodic (as well as quasi-periodic) solutions for \equ(1.1)
is well known, and the fact that no linear term as $\m u$
appears in \equ(1.1) in unessential;
see for instance Refs.~\cita{KP} and \cita{Bo1}, and, very recently,
Ref.~\cita{GY}, where more general nonlinearities are also considered.
Anyway, just because the cases $\m=0$ and $\m \neq 0$ are dealt with
in the same way, all the periodic solutions for $\m=0$ are obtained
as continuations of oscillations involving only one single mode.
Here we consider directly the case $\m=0$, and first we show how
to recover the known results with a different technique,
based on the Lindstedt series method introduced in
Refs.~\cita{GM} and \cita{GMP}, hence we discuss
how to obtain other more complicated periodic solutions
which arise from superposition of several unperturbed modes.

\*

If $\f=0$ every solution of \equ(1.1) can be written as
$$ u(t,x) = \sum_{n=1}^{\io} U_{n} e^{ i n^{2} t} \sin n x =
\sum_{n \in \zzz_{*}} a_{n} e^{ i n^{2} t} e^{i n x}, \qquad a_{-n}=-a_{n} ,
\Eq(1.2) $$
where we have set $\ZZZ_{*}=\ZZZ\setminus\{0\}$.
For $\e\F>0$ we rescale $u \to \sqrt{\e/\F} u$ in \equ(1.1), so obtaining
$$ \cases{
 i u_{t} + u_{xx} = \e |u|^{2}u + O(\e^{2}) , & \cr
u(0,t) = u(\p,t) = 0 , & \cr}
\Eq(1.3)$$
where $O(\e^{2})$ denotes an analytic function of $u$ and $\e$ of
order at least $2$ in $\e$, and we define $\o_{\e}=1+\e$.

We shall consider $\e$ small and we shall show
that for all $m_{0}\in\NNN$ there exists a solution of \equ(1.3),
which is $2\pi/\o_{\e}$-periodic in $t$ and $\e$-close to the function
$$ u_{0}(\o_{\e} t,x) = a(\o_{\e} t, x) - a(\o_{\e} t,-x) ,
\Eq(1.4) $$
with
$$ a(t,x) = A e^{i m_{0}^{2} t} e^{im_{0}x} , \qquad
A = {m_{0} \over \sqrt{3}} ,
\Eq(1.5) $$
provided that $\e$ is in an appropriate Cantor set
(depending on $m_{0}$).

\*

\0{\bf Theorem 1.} {\it Consider the equation \equ(1.1),
where $\f(x)=\F x + O(x^{2})$ is an analytic function,
with $\F \neq 0$. For all $m_{0}\in \NNN$,
define $u_{0}(t,x)=a(t,x)-a(t,-x)$, with $a(t,x)$ as in \equ(1.5).
There are a positive constant $\e_{0}$ and a set $\EE \in[0,\e_{0}]$,
both depending on $m_{0}$, satisfying
$$ \lim_{\e\to0} { {\rm meas}(\EE \cap [0,\e] ) \over \e} = 1 ,
\Eq(1.6) $$
such that for all $\e\in\EE$, by setting $\o_{\e}=1+\e$ and
$$ \|f(t,x)\|_r= \sum_{(n,m)\in \zzz^2}f_{n,m}e^{r(|n|+|m|)}
\Eq(1.7) $$
for analytic $2\pi$-periodic functions,
there exists a $2\p/\o_{\e}$-periodic solution
$u_{\e}(t,x)$ of \equ(1.1), analytic in $(t,x)$, with
$$ \left\| u_{\e}(x,t) - \sqrt{\e/\F} \, u_{0}(\o_{\e}t,x)
\right\|_{\k} \le C \, \e\sqrt{\e} ,
\Eq(1.8) $$
for some constants $C,\k>0$.}

\*

We start by performing a Lyapunov-Schmidt type decomposition:
we look for a solution of \equ(1.3) of the form
$$ \eqalign{
u(x,t) & = \sum_{(n,m) \in \zzz^{2}}
e^{i n \o t + i m x} u_{n,m}
= v(x,t) + w(x,t) , \cr
v(x,t) & = a(\o t,x) - a(\o t,-x) + V(\o t,x) , \cr
V(t,x) & = \sum_{m\in  \zzz} 
e^{i m^{2} t + i m x} V_{m} , \cr
w(x,t) & = \sum_{\scriptstyle (n,m) \in  \zzz^{2} \atop
\scriptstyle n \neq m^{2}}
e^{in \o t +i m x} w_{n,m} , \cr}
\Eq(1.9) $$
with $u_{n,m}\in\RRR$ and $\o=\o_{\e}=1+\e$,
such that one has $w(x,t)=0$ and $V(t,x)=0$ for $\e=0$.
Of course by the symmetry of \equ(1.1), hence of \equ(1.4),
we can look for solutions (if any) which verify
$$ u_{n,m} = - u_{n,-m} , 
\Eq(1.10) $$
for all $n,m\in\ZZZ$.

Inserting \equ(1.9) into \equ(1.3) gives two sets of equations, called
the Q and P equations \cita{CW1}, which are given, respectively, by
$$ \eqalign{
\hbox{Q} \qquad \qquad & 
m^{2} v_{m} = \left[ \f(|v+w|^{2})(v+w) \right]_{m} , \cr
\hbox{P} \qquad \qquad &
\left( \o n - m^2 \right) w_{n,m} = \e \left[
\f(|v+w|^{2})(v+w) \right]_{n,m} , \qquad n \neq m^{2} , \cr}
\Eq(1.11) $$
where we denote by $[F]_{n,m}$ the Fourier component of the function
$F(t,x)$ with labels $(n,m)$, so that
$$ F(x,t) = \sum_{(n,m)\in\zzz^{2}} e^{i n\o t + m x} [F]_{n,m} ,
\Eq(1.12) $$
and we shorthand $[F]_{m}$ the Fourier component of
the function $F(t,x)$ with label $n=m^{2}$; hence
$[F]_{m}=[F]_{m^{2},m}$.

As in Ref.~\cita{GMP} we start by considering the case $\f(x)=x$,
which contains all the relevant features of the problem.
We shall show in Section \secc(4) how to extend the analysis
to more general nonlinearities, which is trivial,
and to more general periodic solutions, which requires some discussion.
The result we obtain at the end is the following one.

\*

\0{\bf Theorem 2.} {\it Consider the equation \equ(1.1),
where $\f(x)=\F x + O(x^{2})$ is an analytic function,
with $\F \neq 0$. For all $N \ge 2$ there are sets of
$N$ positive integers $\MM_{+}$ and sets
of real amplitudes $\{a_{m}\}_{m\in\MM_{+}}$,
such that the following holds. Define
$$ a(t,x) = \sum_{m\in\MM_{+}} e^{im^{2}t+imx}a_{m} ,
\Eq(1.13) $$
and set $u_{0}(x,t)=a(t,x)-a(t,-x)$.
There are a positive constant $\e_{0}$ and a set $\EE \in[0,\e_{0}]$,
both depending on the set $\MM_{+}$, satisfying
$$ \lim_{\e\to0} { {\rm meas}(\EE \cap [0,\e] ) \over \e} = 1 ,
\Eq(1.14) $$
such that for all $\e\in\EE$, by setting $\o_{\e}=1+\e$,
there exist a $2\p/\o_{\e}$-periodic solution
$u_{\e}(t,x)$ of \equ(1.1), analytic in $(t,x)$, with
$$ \left\| u_{\e}(x,t) - \sqrt{\e/\F} \, u_{0}(\o_{\e}t,x)
\right\|_{\k} \le C \, \e\sqrt{\e} ,
\Eq(1.15) $$
for some constants $C,\k>0$.}

\*

In the proof of Theorem 2 a characterization of the sets $\MM_{+}$
will be provided. Hence the proof is constructive. 
What is found is that the integers in $\MM_{+}$ 
have to be large enough and close enough to each other,
so that the solutions which can be continued appear as
wave packets with large Fourier label (wave number).

From a technical point of view the discussion below heavily
relies on \cita{GM} and \cita{GMP}. We confine ourselves to
explain how the renormalization group analysis developed in those
papers applies to the nonlinear Schr\"odinger equation,
by outlining the differences everywhere they appear
and showing how they can be faced.
Hence a full acquaintance with those paper is assumed.
The discussion of Theorem 2 requires some new ideas,
and involves problems which can be considered as
typical of number theory and matrix algebra.

%%%%%%%%%%%%%%%%%%%%%%%%%%%%%%%%%%%%%%%%%%%%%%%%%%%%%%%%%%%%%%%%%%%%%%%%%%
%%%%%%%%%%%%%%%%%%%%%%%%%%%%%%%%%%%%%%%%%%%%%%%%%%%%%%%%%%%%%%%%%%%%%%%%%%
\*\*
%%%%%%%%%%%%%%%%%%%%%%%%%%%%%%%%%%%%%%%%%%%%%%%%%%%%%%%%%%%%%%%%%%%%%%%%%%
%%%%%%%%%%%%%%%%%%%%%%%%%%%%%%%%%%%%%%%%%%%%%%%%%%%%%%%%%%%%%%%%%%%%%%%%%%
\section(2,Lindstedt series expansion)
%%%%%%%%%%%%%%%%%%%%%%%%%%%%%%%%%%%%%%%%%%%%%%%%%%%%%%%%%%%%%%%%%%%%%%%%%%
%%%%%%%%%%%%%%%%%%%%%%%%%%%%%%%%%%%%%%%%%%%%%%%%%%%%%%%%%%%%%%%%%%%%%%%%%%

\0Given a sequence
$\{\n_{m}(\e)\}_{|m|\ge 1}$, such that $\n_{m}=\n_{-m} $,
we define the {\rm renormalized frequencies} as
$$ \tilde\o_m^2 \= \o_m^2 - \n_m , \qquad \o_{m}=|m| ,
\Eq(2.1) $$
and the quantities $\n_{m}$ will be called the counterterms.

\*

By the above definition and the parity properties \equ(1.10)
the P equation in \equ(1.11) can be rewritten as
$$ \eqalign{
\left( \o n - \tilde\o_m^2 \right) w_{n,m} 
& = \n_{m} w_{n,m} + \e [\f(v+w)]_{n,m} \cr
& = \n_{m}^{(a)} w_{n,m} + \n_{m}^{(b)} w_{n,-m} + \e [\f(v+w)]_{n,m}, \cr}
\Eq(2.2) $$
where
$$ \n_{m}^{(a)}-\n_{m}^{(b)}=\n_{m} .
\Eq(2.3) $$

With the notations of \equ(1.12), and recalling that we are
considering $\f(x)=x$, we have
$$ \eqalign{
\left[ |v+w|^{2}(v+w) \right]_{m}
& = [|v|^{2}v]_{m}+[|w|^{2}w]_{m}+
[2|v|^{2} w + \lis w v^{2}]_{m} +
[2|w|^{2} v + \lis v w^{2}]_{m} \cr
& \equiv [|v|^{2}v]_{m} + [G_{2}(v,w)]_{m} , \cr}
\Eq(2.4) $$
where $G_{2}(v,w)$ is at least linear in $w$.

We can write $v=a+b+V$, with $b(\o t,x)=-a(\o t,-x)$, so that
$$ b(t,x) = B e^{i m_{0}^{2} t - i m_{0}x} , \qquad B = - A ,
\Eq(2.5) $$
so that we obtain
$$ \eqalign{
[|v|^{2}v]_{m} & =
|A| ^{2} A \d_{m,m_{0}} + |B| ^{2} B \d_{m,-m_{0}} +
2 |A|^{2}B \d_{m,-m_{0}} + 2 |B|^{2} A \d_{m,m_{0}} \cr
& + 2 |A|^{2} V_{m} + 2|B|^{2} V_{m} + 2 \lis A B V_{-m} \d_{m,-m_{0}}
+ 2 \lis B A V_{-m} \d_{m,m_{0}} \cr
& + \lis V_{m} A^{2} \d_{m,m_{0}} + 2 \lis V_{-m} AB \d_{m,\pm m_{0}}
+ \lis V_{m} B^{2} \d_{m,-m_{0}} + [G_{1}(v)]_{m} , \cr}
\Eq(2.6) $$
where $G_{1}(v)$ is at least quadratic in $V$.

Then, by setting $G(v,w)=G_{1}(v)+G_{2}(v,w)$,
the Q equation in \equ(1.11) can be rewritten for $m=m_{0}$ as
$$ \eqalign{
m_{0}^{2} A & = |A|^{2} A + 2 |B|^{2}A , \cr
m_{0}^{2} V_{m_{0}} & =
2 |A|^{2} V_{m_{0}} + 2|B|^{2} V_{m_{0}} + \lis V_{m_{0}} A^{2} +
2 AB \lis V_{-m_{0}} + 2 \lis B A V_{-m_{0}} +
[G(v,w)]_{m_{0}} , \cr}
\Eq(2.7) $$
and for positive $m \neq m_{0}$ as
$$ m^{2} V_{m} = 2 |A|^{2} V_{m} + 2|B|^{2} V_{m} + [G(v,w)]_{m} ,
\Eq(2.8) $$
while the equation for negative values of $m$ can be obtained
by using the symmetry properties \equ(1.10),
which imply $V_{-m}=-V_{m}$.

By defining $\a=|A|^{2} = \lis A A$ and using the identities
$$ \a = \lis A A = \lis B B = - \lis A B = -\lis B A , \qquad
\b = A A = B B = - A B = - B A ,
\Eq(2.9) $$
which follow trivially from the definitions of $A$ and $B$
in \equ(1.5) and \equ(2.5) respectively, we can rewrite \equ(2.7) as
$$ \eqalign{
m_{0}^{2} A & = 3 \a A , \cr
m_{0}^{2} V_{m_{0}} & =
4 \a V_{m_{0}} + \b \lis V_{m_{0}} -
2 \b \lis V_{-m_{0}} - 2 \a V_{-m_{0}} +
[G(v,w)]_{m_{0}} , \cr}
\Eq(2.10) $$
where the first equation defines $\a=|A|^{2} = m_{0}^{2}/3$.
By using once more the identities \equ(1.10) and imposing that
the coefficients $V_{m}$ be real, so that
$\a=\b$ in \equ(2.9), we can write
the second equation in \equ(2.10) and the equation \equ(2.8),
respectively, as
$$ \cases{
m_{0}^{2} V_{m_{0}} = 9 \a V_{m_{0}} + [G(v,w)]_{m_{0}} , & \cr
m^{2} V_{m} = 4\a V_{m} + [G(v,w)]_{m} , & \cr}
\Eq(2.11) $$
so that we find
$$ \cases{
V_{m_{0}} = {\displaystyle - {1 \over 6\a} [G(v,w)]_{m_{0}} } , & \cr
& \cr
V_{m} = {\displaystyle -{1 \over \a} [G(v,w)]_{m} } , & \cr}
\Eq(2.12) $$
respectively for positive $m_{0}$ and $m\neq m_{0}$.

Finally we write
$$ w_{n,m} = g(n,m) \left(
\m \n_{m}^{(a)} w_{n,m} + \m \n_{m}^{(b)} w_{n,-m} +
\m \e [|v+w|^{2}(v+w)]_{n,m} \right) ,
\Eq(2.13) $$
where
$$ g(n,m) = { 1 \over \o n - \tilde \o_{m}^{2} } ,
\qquad n \neq m^{2} ,
\Eq(2.14) $$
and we look for a solution $u_{n,m}$ in the form of a 
power series expansion in $\m$,
$$ u_{n,m} = \sum_{k=0}^{\io} \m^{k} u_{n,m}^{(k)} ,
\Eq(2.15) $$
with $u_{n,m}^{(k)}$ depending on $\e$ and on
the parameters $\n_{m}^{(c)}$, with $c=a,b$ and $|m|\ge 1$.

So we obtain recursive definitions of the coefficients
$u^{(k)}_{n,m}$. The coefficients $w^{(k)}_{n,m}$
verify for $k\ge1$ the recursive equations
$$ w^{(k)}_{n,m} = g(n,m) \left(
\n_{m}^{(a)} w_{n,m}^{(k-1)} + \n_{m}^{(b)} w_{n,-m}^{(k-1)} +
[|v+w|^{2}(v+w)]^{(k-1)}_{n,m} \right) ,
\Eq(2.16) $$
where
$$ [|v+w|^{2}(v+w)]^{(k)}_{n,m}  = \sum_{k_{1}+k_{2}+k_{3}=k}
\sum_{\scriptstyle - n_{1}+n_{2}+n_{3} = n \atop
\scriptstyle - m_{1}+m_{2}+m_{3} = m} \lis u^{(k_{1})}_{n_{1},m_{1}}
u^{(k_{2})}_{n_{2},m_{2}} u^{(k_{3})}_{n_{3},m_{3}} ,
\Eq(2.17) $$
with
$$ u^{(0)}_{n,m} = \cases{
A , & if $n=m^{2}$ and $m>0$, \cr
B , & if $n=m^{2}$ and $m<0$, \cr
0 , & otherwise \cr}
\Eq(2.18) $$
and, for $k\ge 1$,
$$ u^{(k)}_{n,m} = \cases{
V^{(k)}_{m} , & if $n=m^{2}$ , \cr
w^{(k)}_{n,m} , & if $n \neq m^{2}$ , \cr}
\Eq(2.19) $$
while the coefficients $V_{m}^{(k)}$ verify for $k\ge 1$ the equations
$$ V^{(k)}_{m} = g(m^{2},m) =
\sum_{k_{1}+k_{2}+k_{3}=k}
\mathop{{\sum}^{*}}_{\scriptstyle - n_{1}+n_{2}+n_{3} = m \atop
\scriptstyle - m_{1}+m_{2}+m_{3} = m} \lis u^{(k_{1})}_{n_{1},m_{1}}
u^{(k_{2})}_{n_{2},m_{2}} u^{(k_{3})}_{n_{3},m_{3}} ,
\Eq(2.20) $$
where
$$ g(m^{2},m)=\cases{
- {\displaystyle {1 \over 18 m_{0}^{2}}} , & if $|m|=m_{0}$ , \cr
- {\displaystyle {1 \over 3 m_{0}^{2}}} , & if $|m| \neq m_{0}$ , \cr}
\Eq(2.21) $$
and the $*$ means that there appear only contributions
either with at least one coefficient with $n\neq m^{2}$ 
or with at least two labels $k_{i}\geq 1$.

To prove Theorem 1 we can proceed in two steps as in Ref.~\cita{GMP}.
The first step consists in looking for the solution
of the recursive equations
by considering $\tilde\o=\{\tilde\o_{m}\}_{|m|\ge 1}$
as a given set of parameters satisfying the Diophantine conditions
(called respectively the first and the second Mel$^{\prime}$nikov conditions)
$$ \eqalignno{
& \left| \o n \pm \tilde \o_{m}^{2} \right|
\ge C_{0} |n|^{-\t}
\qquad \forall n\in\ZZZ_{*} \hbox{ and }
\forall m \in \ZZZ_{*}
\hbox{ such that } n \neq m^{2} ,
& \eq(2.22) \cr
& \left| \o  n \pm \left( \tilde\o_{m}^{2} \pm \tilde\o_{m'}^{2}
\right) \right| \ge C_{0} |n|^{-\t} \qquad
\forall n\in\ZZZ_{*} \hbox{ and }
\forall m,m' \in \ZZZ_{*}
\hbox{ such that } |n|\neq |m^{2} \pm (m')^{2}| , \cr} $$
with positive constants $C_{0},\t$. We can assume without loss
of generality $C_{0}\le 1/2$.
%If $\t \le 2$ then one can require only the first Mel'nikov
%conditions in \equ(2.22), as in Ref.~\cita{GMP}.

We shall show in Section \secc(3) how to adapt the discussion
in Ref.~\cita{GMP} in order to obtain the following result.

\*

\0{\bf Proposition 1.} {\it Consider a sequence
$\tilde\o=\{\tilde\o_m\}_{|m| \ge 1}$ verifying \equ(2.22),
with $\o=\o_{\e}=1+\e$ and such that
$|\tilde\o_{m}^{2}-m^{2}|\le C_{1}\e$ for some constant $C_{1}$.
For all $\m_{0}>0$ there exists $\e_{0}>0$ such that
for $|\m|\le \m_{0}$ and $0<\e<\e_{0}$ there is
a sequence $\n(\tilde\o,\e;\m) = \{\n_{m}(\tilde\o,\e;\m)\}_{|m| \ge 1}$,
where each $\n_{m}(\tilde\o,\e;\m)$ is analytic in $\m$,
such that there are coefficients $u^{(k)}_{n,m}$
which solve the recursive equations \equ(2.16) and \equ(2.19),
with $\n_{m}=\n_{m}(\tilde\o,\e)$,
and define a function $u(t,x;\tilde\o,\e;\m)$ which is analytic in $\m$,
analytic in $(t,x)$ and $2\pi/\o_{\e}$-periodic in $t$.}

\*

Then in Proposition 1 one can fix $\m_{0}=1$, so that
one can choose $\m=1$ and set
$u(t,x;\tilde\o,\e)=u(t,x;\tilde\o,\e;1)$ and
$\n_{m}(\tilde\o,\e)=\n_{m}(\tilde\o,\e;1)$.

\*

The second step, also to be proved in Section \secc(3),
consists in inverting \equ(2.1), with $\n_{m}=\n_{m}(\tilde\o,\e)$
and $\tilde\o$ verifying \equ(2.22). This requires
some preliminary conditions on $\e$, given by
the Diophantine conditions
$$ \left| \o n \pm m \right|
\ge c \, C_{0} |n|^{-\t_{0}}
\qquad \forall n\in\ZZZ_{*} \hbox{ and }
\forall m \in \ZZZ_{*} 
\hbox{ such that } n \neq m ,
\Eq(2.23) $$
with positive constants $c>1$ and $\t_{0}>1$.
Then we can solve iteratively \equ(2.1),
by imposing further non-resonance conditions besides \equ(2.23).
At each iterative step one has to exclude
some further values of $\e$, and at the end the left values 
fill a Cantor set $\EE$ with large relative measure in $[0,\e_{0}]$
and $\tilde\o$ verify \equ(2.22).

The result of this second step can be summarized as follows.

\*

\0{\bf Proposition 2.} {\it
There are $\d>0$ and a set $\EE\subset [0,\e_{0}]$
with complement of relative Lebesgue measure of order $\e_{0}^{\d}$
such that for all $\e\in\EE$ there exists
$\tilde\o=\tilde\o(\e)$ which solves \equ(2.1)
and satisfy the Diophantine conditions \equ(2.22) with
$|\tilde\o_{m}^{2}-m^{2}|\le C_{1}\e$ for some positive constant $C_{1}$.}

\*

The proof follows the same strategy as in Ref.~\cita{GMP},
which we refer to for further details. The slight changes
will be briefly discussed in Section \secc(3).

%%%%%%%%%%%%%%%%%%%%%%%%%%%%%%%%%%%%%%%%%%%%%%%%%%%%%%%%%%%%%%%%%%%%%%%%%%
%%%%%%%%%%%%%%%%%%%%%%%%%%%%%%%%%%%%%%%%%%%%%%%%%%%%%%%%%%%%%%%%%%%%%%%%%%
\*\*
%%%%%%%%%%%%%%%%%%%%%%%%%%%%%%%%%%%%%%%%%%%%%%%%%%%%%%%%%%%%%%%%%%%%%%%%%%
%%%%%%%%%%%%%%%%%%%%%%%%%%%%%%%%%%%%%%%%%%%%%%%%%%%%%%%%%%%%%%%%%%%%%%%%%%
\section(3,Renormalization and proof of Theorem 1)
%%%%%%%%%%%%%%%%%%%%%%%%%%%%%%%%%%%%%%%%%%%%%%%%%%%%%%%%%%%%%%%%%%%%%%%%%%
%%%%%%%%%%%%%%%%%%%%%%%%%%%%%%%%%%%%%%%%%%%%%%%%%%%%%%%%%%%%%%%%%%%%%%%%%%

\0We refer to Section 3 in Ref.~\cita{GMP} for the tree expansion.
With respect to that paper the following changes have to be
performed. Each line $\ell$ carries a momentum
$(n_{\ell},m_{\ell})$ and label $\g_{\ell}$,
with $\g_{\ell}=v$ if $n_{\ell}=m_{\ell}^{2}$
and $\g_{\ell}=w$ otherwise. The corresponding propagator is
given by $g_{\ell}=g(n_{\ell},m_{\ell})$, with
$g(n,m)$ defined in \equ(2.14) and \equ(2.21)
if the line comes out from a node, while it is $g_{\ell}=1$
if it comes out from an end-point.

If we denote by $s_{\vvvv}$ the number of lines
entering the node $\vvv$ one can have either $s_{\vvvv}=1$
or $s_{\vvvv}=3$. In the latter we call $L_{\vvvv}$ the set
of lines entering $\vvv$: we associate to each line
$\ell\in L_{\vvvv}$ a label $s(\ell)\in\{\pm1\}$ with
the constraint $\sum_{\ell\in L_{\vvvv}}s(\ell)=1$.

To each end-point $\vvv$ a mode label $(n_{\vvvv},m_{\vvvv})$
is associated, with $m_{\vvvv}=\pm m_{0}$ and $n_{\vvvv}=m_{0}^{2}$.

The momentum $(n_{\ell},m_{\ell})$ of a line $\ell=\ell_{\vvvv}$
coming out from a node $\vvv$ is given by
$$ n_{\ell} = \sum_{\wwww \in E(\th) \atop \wwww \preceq \vvvv}
(-1)^{S(\wwww,\ell)} n_{\wwww} ,
\qquad m_{\ell} = \sum_{\wwww \in E(\th) \atop \wwww \preceq \vvvv}
(-1)^{S(\wwww,\ell)} m_{\wwww} +
\sum_{\wwww \in V_{w}^{1}(\th) : c_{\wwwww}=b \atop \wwww \preceq \vvvv}
(-2m_{\ell_{\wwww}}) ,
\Eq(3.1) $$
where $S(\wwww,\ell)$ is the number of lines $\ell$ with $s(\ell)=-1$
between $\www$ and $\ell$,
and the node factor $\h_{\vvvv}$ is defined as
$$ \h_{\vvvv} = \cases{
1/3 , & $\vvv \in V_{v}(\th)$ , \cr
\e , & $\vvv \in V_{w}^{3}(\th)$ , \cr
\n_{m_{\ell_{\vvvv}}}^{(c_{\vvvv})} , & $\vvv \in V_{w}^{1}(\th)$ , \cr}
\Eq(3.2) $$
and $s_{\vvvv}=3$ for all $\vvv\in V_{v}$ (so that $V_{v}^{1}=\emptyset$
in the present case). Finally the end-point factor is
$$ V_{\vvvv} = \cases{
A , & $m_{\vvvv} = m_{0}$ , \cr
-A , & $m_{\vvvv} = -m_{0}$ . \cr}
\Eq(3.3) $$
All the other notations used below are as in Ref.~\cita{GMP}.

Introducing a multiscale decomposition as in Section 4 of Ref.~\cita{GMP}
we can define for the lines $\ell$ with $\g_{\ell}=w$
the propagator on scale $h\ge -1$ as
$$ g_{\ell}^{(h)} = \chi_{h}(|\o n_{\ell}-\tilde \o_{m}^{2}|) \,
g_{\ell} = { \chi_{h}(|\o n_{\ell}-\tilde \o_{m}^{2}|)
\over \o n_{\ell} - \tilde \o_{m}^{2} } ,
\Eq(3.4) $$
where $\chi_{h}(x)$ is a $C^{\io}$ function nonvanishing
for $2^{-h-1}C_{0}<|x|<2^{-h+1}C_{0}$ if $h\ge 0$
and for $|x|>C_{0}$ if $h=-1$.

Then for each tree $\th$ one can define the tree value as
$$ \Val(\th) =
\left( \prod_{\ell\in L(\th)} g_{\ell}^{(h_{\ell})} \right)
\left( \prod_{\vvvv\in V(\th)} \h_{\vvvv} \right)
\left( \prod_{\vvvv\in E(\th)} V_{\vvvv} \right) ,
\Eq(3.5) $$
so that one has
$$ u^{(k)}_{n,m} = \sum_{\th \in \Theta^{(k)}_{n,m}} \Val(\th) ,
\Eq(3.6) $$
where $\Th^{(k)}_{n,m}$ is the set of tress $\th$ of order $k$,
that is with $|V_{w}(\th)|=k$, and with momentum $(n,m)$ associated
to the root line. Note that one has $|V_{v}(\th)| \le 2|V_{w}(\th)|=2k$
and $|E(\th)|\le 2(|V_{w}(\th)|+|V_{v}(\th)|)+1 \le 6k+1$
(see Lemma 3 in Ref.~\cita{GMP}).

Clusters and self-energy graphs are defined as in Ref.~\cita{GMP}.
In particular the self-energy value is given by
$$ \VV_{T}^{h}(\o n,m)=
\Big( \prod_{\ell \in T} g^{(h_{\ell})}_{\ell} \Big)
\Big( \prod_{\vvvv \in V(T)} \h_{\vvvv} \Big)
\Big( \prod_{\vvv \in E(T) } V_{\vvvv} \Big) ,
\Eq(3.7) $$
where $h=h^{(e)}_T$ is the minimum between the scales of the
two external lines of $T$ (they can differ at most by a unit
and $h^{(e)}_{T}\ge 0$), and, given a self-energy graph, one has
$$ \eqalign{
n(T) & \= \sum_{\vvvv \in E(T)}
(-1)^{S(\vvvv,\ell_{T}^{1})} n_{\vvvv} = 0 , \cr
m(T) & \= \sum_{\vvvv \in E(T)}
(-1)^{S(\vvvv,\ell_{T}^{1})} m_{\wwww} + m_{\vvvv} +
\sum_{\scriptstyle \wwww \in V_{w}^{1}(T) \atop
\scriptstyle c_{\wwww}=b} \left( -
2m_{\ell_{\wwww}} \right) \in \{0,2m\} , \cr}
\Eq(3.8) $$
by definition of self-energy graph (recall that $\ell_{T}^{1}$
is the exiting line of $T$). One says that $T$ is
a self-energy graph of type $c=a$ when $m(T)=0$ and a resonance
of type $c=b$ when $m(T)=2m$. 

The following results hold.

\*

\0{\bf Lemma 1.} {\it Assume that there is a constant $C_{1}$
such that $|\tilde\o_{m}^{2}-m^{2}|< C_{1}\e$ for all $m\ge1$.
If $|\o n_{\ell} - \tilde \o_{m}^{2} | < 1/2$ and $\e$
is small enough then $\min\{ n_{\ell},m_{\ell}^{2}\} > 1/4\e$.}

\*

\0{\it Proof.} One has $\o n - \tilde \o_{m}^{2} =
\e n + (n - m^{2}) + \n_{m}$, so that
$|\o n - \tilde \o_{m}^{2}|>1/2$ for $n\neq m^{2}$
and $0<n<1/3\e$. Moreover if $|\o n - \tilde \o_{m}^{2}|<1/2$
then one has $n>0$ and $m^{2}> \o n - |\n_{m}| - 1/2 > 1/4\e$. \qed

\*

Hence if $n_{\ell} < 1/4\e$ we can bound $|g(n_{\ell},m_{\ell})|\le 2$
while if $n_{\ell} \ge 1/4\e$ in general we can bound
$|g(n_{\ell},m_{\ell})| \le 2^{h+1}C_{0}^{-1}$.
To any line $\ell$ with $n_{\ell}<1/4\e$ we can assign
a scale label $h_{\ell}=-1$.

\*

\0{\bf Lemma 2.} {\it Assume that there is a constant $C_{1}$
such that $|\tilde\o_{m}^{2}-m^{2}|< C_{1}\e$.
Define $h_{0}$ such that $2^{h_{0}} < 16 C_{0} /\sqrt{\e} < 2^{h_{0}+1}$.
Then for $h \ge h_{0}$ one has
$$ \lis N_{h}(\th) \le 4 k 2^{(2-h)/\t} -C_{h}(\th) + S_{h}(\th)
+M^{\n}_{h}(\th) ,
\Eq(3.9) $$
with the same notations as in Ref.~{\rm \cita{GMP}}.}

\*

\0{\it Proof.} The proof as for Lemma 5 of Ref.~\cita{GMP}.
Again the only case which deserves attention is when
one has a cluster $T$ with two external lines $\ell$ and
$\ell_{1}$ both on scales $\ge h$, so that, with the
same notations as in Ref.~\cita{GMP}, one has
$$ 2^{-h+2} C_{0} \ge \big| \o (n_{\ell}-n_{\ell_{1}}) +
\h_{\ell} \tilde\o_{m_{\ell}}^{2} +
\h_{\ell_{1}} \tilde\o_{m_{\ell_{1}}}^{2} \big| .
\Eq(3.10) $$
Then $|n_{\ell}-n_{\ell_{1}}|=|m_{\ell}^{2} \pm m_{\ell_{1}}^{2}|$
would require $|n_{\ell}-n_{\ell_{1}}| \geq |m_{\ell}|+
|m_{\ell_{1}}| > 1/\sqrt{\e}$, while \equ(3.10) would become
$2^{-h+2}C_{0}>|\e(n_{\ell}-n_{\ell_{1}})|-2C_{1}\e$.
Combining the two inequalities one would obtain $C_{0}2^{-h+3} >
\sqrt{\e}$, which contradicts the condition $h \ge h_{0}$.
Then one proceeds as in Ref.~\cita{GMP}.\qed

\*

\0{\bf Lemma 3.} {\it Assume that there is a constant $C_{1}$
such that $|\tilde\o_{m}^{2}-m^{2}|< C_{1}\e$.
Then one has
$$ \prod_{h=0}^{h_{0}-1}
\prod_{\ell\in L(\th) \atop h_{\ell}=h} |g^{(h_{\ell})}_{\ell}|
\le C_{2}^{k} \e^{-k/2} ,
\Eq(3.11) $$
for some positive constant $C_{2}$.}

\*

\0{\it Proof.} If $h<h_{0}$ one has $|g^{(h_{\ell})}_{\ell}|
\le C_{0}2^{-h_{0}+1} < \sqrt{\e}/4$, and the number of lines
$\ell$ with scales $0 \ge h_{\ell}<h_{0}$ can be bounded
by the total number of lines $\ell$ with label $\g_{\ell}=w$,
which is less than $k$. \qed

\*

The renormalized expansion is defined as in Section 5 Ref.~\cita{GMP},
with the only difference that now the action of the localization
operator $\LL$ is such that
$$ \LL \VV_{T}^{h}(\o n, m) = \VV_{T}^{h}(\tilde\o_{m}^{2},m) ,
\Eq(3.12) $$
so that, in the definition of the set $E_{0}(\th)$
(see item (7$^{\prime}$) in Section 5 of Ref.\cita{GMP}),
we set $\lis\o_{m}=\tilde\o_{m}^{2}/\o$. Up to these
notational changes no real difference appears with respect
to the discussion carried out in Ref.\cita{GMP}.
Therefore by using the lemmata above, the discussion proceeds as in
Ref.~\cita{GMP}, hence we omit the details of the proof of Proposition 1.

In order to apply the results stated above we have to prove
the following result (as discussed in Ref.~\cita{GMP}).

\*

\0{\bf Lemma 4.} {\it There exists a positive constant $C_{3}$ such that
one has $|\n_{m}(\tilde\o,\e)| < C_{3}\e$.}

\*

\0{\it Proof.} As in Lemma 16 in \cita{GMP}.

\*

The construction of the perturbed frequencies $\tilde\o_{m}$
and the bound of the measure of the admissible values of $\e$,
in order to prove Proposition 2,
proceeds as in \cita{GMP}, with some minor differences
(which are in fact simplifications) that we outline below.

The condition \equ(2.23) on $\e$ can be imposed exactly
as in Ref.~\cita{GMP}, and $c$ can be chosen as $c=2$.

To impose the Mel'nikov conditions in \equ(2.22)
one has to use that $|\dpr f/\dpr\e| \ge n/2$,
if the function $f(\e)$ is defined through
$$ f(\e(t)) \= (1+\e(t)) n - \tilde\o_{m}^{2}(\e(t)) =
t {C_{0} \over |n|^{\t}} , \qquad t \in [-1,1] .
\Eq(3.13) $$
when dealing with the first Mel'nikov conditions, and through
$$ f(\e(t)) \= (1+\e(t)) |n_{\ell}|-|\tilde\o_{m}^{2}(\e(t)) \pm
\tilde \o_{m'}^{2}(\e(t))| =
t {C_{0} \over |n|^{\t}} , \qquad t \in [-1,1] .
\Eq(3.14) $$
when dealing with the second Mel'nikov conditions.

In the case of the first Mel'nikov conditions one has
to consider only the values of $n$ such that
$n \ge N_{0}=O(\e_{0}^{-1/\t_{0}})$, as $|\n_{m}| < C_{1}\e$,
and for each $n$ the set $\MM_{0}(n)$ of $m$'s such that the
conditions are not satisfied contains at most $2+\e\sqrt{n}$ values.
Therefore the measure of the set of excluded values of $\e$
turns out to be bounded proportionally to
$$ \sum_{n> N_{0}} {C_{0} \over n^{\t+1}} \left( 2 + \e\sqrt{n} \right)
\le \hbox{const.} \, \e_{0}^{1+\d_{1}} ,
\Eq(3.15) $$
with $\d_{1}>0$ provided that one takes $\t>\t_{0}$.

In the case of the second Mel'nikov conditions one has
to use that if $n$ is close to $|m^{2} - (m')^{2}|$
then $|n|$ is of order $\left||m|-|m'|\right|(|m|+|m'|)$, with
$|m|-|m'|\neq0$, so that $|m|+|m'| \le n$.
This means that for each $n$ the number of pairs
$(m,m')$ one has to sum over is at most proportional to $|n|$.
The same happens (trivially) when $n$ is close
to $m^{2}+(m')^{2}$. In both cases one has to sum only
on the values of $n$ such that $n \ge N_{0}=O(\e_{0}^{-1/\t_{0}})$,
so that one have to exclude a set of values of $\e$
whose measure is bounded proportionally to
$$ \sum_{n> N_{0}} {C_{0} \over n^{\t+1}} n 
\le \hbox{const.} \, \e_{0}^{1+\d_{2}} ,
\Eq(3.16) $$
with $\d_{2}>0$ provided that one takes $\t>\t_{0}+1$.

At the end one has to choose $\t>\t_{0}+1>2$.
Again we refer to Ref.~\cita{GMP} for further details.

%%%%%%%%%%%%%%%%%%%%%%%%%%%%%%%%%%%%%%%%%%%%%%%%%%%%%%%%%%%%%%%%%%%%%%%%%%
%%%%%%%%%%%%%%%%%%%%%%%%%%%%%%%%%%%%%%%%%%%%%%%%%%%%%%%%%%%%%%%%%%%%%%%%%%
\*\*
%%%%%%%%%%%%%%%%%%%%%%%%%%%%%%%%%%%%%%%%%%%%%%%%%%%%%%%%%%%%%%%%%%%%%%%%%%
%%%%%%%%%%%%%%%%%%%%%%%%%%%%%%%%%%%%%%%%%%%%%%%%%%%%%%%%%%%%%%%%%%%%%%%%%%
\section(4,Extension of the results)
%%%%%%%%%%%%%%%%%%%%%%%%%%%%%%%%%%%%%%%%%%%%%%%%%%%%%%%%%%%%%%%%%%%%%%%%%%
%%%%%%%%%%%%%%%%%%%%%%%%%%%%%%%%%%%%%%%%%%%%%%%%%%%%%%%%%%%%%%%%%%%%%%%%%%

\0The extension of the results of the previous sections to the case
in which $\f(x)$ is any analytic function with $\f'(0)\neq0$,
can be easily dealt with by reasoning as in Ref.~\cita{GMP}.

So we pass directly to discuss the case of more general periodic
solutions to be continued.

For $\e=0$ (and $\f(x)=x$, again for simplicity)
we call $v_{0}=a+b$ the solution of the Q equation
for $\e=0$, by writing
$$ a(t,x) = \sum_{m=1}^{\io} a_{m} e^{im^{2}t+imx} ,
\Eq(4.1) $$
with coefficients $a_{m}$ to be determined, and
setting $b(t,x)=-a(t,-x)$.
Then the Q equation in \equ(1.1) becomes
$$ \eqalign{
m^{2} v_{0,m} & = \sum_{
\scriptstyle -m_{1}+m_{2}+m_{3} = m \atop
\scriptstyle -m_{1}^{2}+m_{2}^{2}+m_{3}^{2} = m^{2} }
\lis v_{0,m_{1}} v_{0,m_{2}} v_{0,m_{3}} \cr
& = 2 v_{0,m} \sum_{ m' \neq m} \lis v_{0,m'} v_{0,m'} +
\lis v_{0,m} v_{0,m} v_{0,m} , \cr}
\Eq(4.2) $$
so that we obtain
$$ v_{0,m} \left( m^{2} - 2\|v_{0}\|^{2} + |v_{0,m}|^{2} \right) = 0 ,
\Eq(4.3) $$
where we have defined
$$ \MM = \left\{ m \in \zzz : v_{0,m} \neq 0 \right\} , \qquad
\MM_{+} = \left\{ m \in \MM : m >0 \right\} ,
\Eq(4.4) $$
and set
$$ \|v_{0}\|^{2} \= \sum_{m\in\zzz} v_{0,m}^{2} =
\sum_{m\in\MM} v_{0,m}^{2} .
\Eq(4.5) $$
Hence \equ(4.3) can be satisfied either if $v_{0,m}=0$ or,
when $v_{0,m} \neq 0$, if
$$ \| v_{0} \|^{2} = {2M \over 4N -1} ,
\Eq(4.6) $$
where we have set
$$ 2N = |\MM| = \#\left\{ m \in \MM \right\} , 
\qquad 2M = \sum_{m\in\MM} m^{2} .
\Eq(4.7) $$
By inserting \equ(4.6) into \equ(4.3), setting
$$ a_{m} = v_{0,m} , \quad m>0 , \qquad \| a \|^{2} =
\sum_{m \in \MM_{+}} a_{m}^{2} = {1\over2} \| v_{0} \|^{2} ,
\Eq(4.8) $$
and writing $\MM_{+} = \{ m_{1}, m_{2}, \ldots, m_{N} \}$,
with $m_{k}<m_{k+1}$, $k=1,\ldots,N-1$, we obtain
$$ a_{m_{k}}^{2} = 4\| a \|^{2} - m_{k}^{2} =
{4\over 4N -1} \left( m_{1}^{2} + m_{2}^{2} + \ldots +
m_{N}^{2} \right) - m_{k}^{2} , \qquad k=1,\ldots,N ,
\Eq(4.9) $$
which makes sense as long as
$$ \max_{m \in \MM_{+}} m^{2} \le {4 \over 4N-1}
\sum_{m\in\MM_{+}} m^{2} .
\Eq(4.10) $$
The following result is easily proved.

\*

\0{\bf Lemma 5.} {\it For all $N \ge 2$ there are solutions of \equ(4.9)
such that $4\|a\|^{2}$ is not an integer.}

\*

\0{\it Proof.} To obtain a solution one can take
$m_{k}=m_{N}-(N-k)$ for $k=1,\ldots,N$,
and choose $m_{N} \ge 4N(N-1)$. Choose $m_{N}=(4N+j)(N-1)$,
with $j\in\{0,1\}$: then $4(m_{1}^{2}+\ldots+m_{N}^{2})$
can not be a multiple of $4N-1$ for both $j=0$ and $j=1$. \qed

\*

Here we are confined ourselves only to an existence result.
Of course more general solutions can be envisaged, with more
spacing between the involved wave numbers $m_{k}$.
The result above can indeed be strenghtened as follows.

\*

\0{\bf Lemma 6.} {\it For all $N \ge 2$ and for all increasing lists
of positive integers $I:=\{i_{1},\ldots,i_{N-1}\}$ there exists $m_N(I)$
($m_N$ for short) such that \equ(4.10) has a solution in the
set $\MM_{+}=\{m_{N}-i_{N-1},\ldots,m_{N}-i_{1},m_{N}\}$
with $4\|a\|^{2} \neq m^{2}$ for all $m\notin\MM$.}

\*

\0{\it Proof.} Fix the set of integers $I=\{i_{1},\ldots,i_{N-1}\}$,
and consider the expression $M-(N-1/4)j^{2}$ for $j\in\NNN$.
For $j=m_{N}$ it becomes a polynomial of degree two in $m_{N}$,
with positive leading coefficient $1/4$ and positive discriminant.
Hence there is an integer $K_{1}$ such that for all $m_{N}>K_{1}$
one has $f_{1}(m_{N})\= M-(N-1/4)m_{N}^{2}>0$,
hence \equ(4.10) is satisfied.

The inequality $M-(N-1/4)j^{2} > 0$ is trivially satisfied
for $j \le m_{N}$, so that it is enough to look
an integer $m_{N}>K_{1}$ such that one has
$f_{2}(m_{N})\=M-(N-1/4)(m_{N}+1)^{2}<0$. Again
$f_{2}$ is a polynomial of degree two in $m_{N}$,
with positive leading coefficient $1/4$
and positive discriminant, so that there
exist two integers $K_{2}< K_{3}$ such that $f_{2}(m_{N})<0$
for $K_{2}< m_{N} < K_{3}$. Moreover
$K_{3}-K_{2}\ge 4(N-1)$, so that
there is $m_{N}$ satisfying \equ(4.10) such that $4\|a\| \neq j^{2}$
for all $j\in\NNN$. \qed

\*

To have solutions of \equ(4.9) requires the integers in $\MM_{+}$ to
be large enough, and not too distant from each other. Hence the
solutions whose existence is stated in Lemma 5 have the form
of wave packets centered around some harmonic large enough.

Hence we have proved the following result.

\*

\0{\bf Lemma 7.} {\it For any $N$ there are sets $\MM$
and functions $v_{0}(x,t)=a(\o t,x)-a(\o t,-x)$, with
$$ a(t,x)  = \sum_{m\in\MM_{+}} e^{i m^{2} \o t +imx} a_{m} ,
\Eq(4.11) $$
which solve the Q equation with $\e=0$.}

\*

Moreover, by using once more the parity properties $V_{-m}=-V_{m}$,
one obtains for $m\in\MM_{+}$
$$ \sum_{m'\in \MM_{+}} \AAA_{m,m'} V_{m'} =\left[ G(v,w)
\right]_{m} ,
\Eq(4.12) $$
where $\AAA$ is an $N\times N$ matrix with entries
$$ \AAA_{m,m'}= \cases{
m^{2} -  4 \|a\|^{2} - 5 a_{m}^{2} , & $m=m'$ , \cr
-8a_{m}a_{m'} , & $m\neq m'$ , \cr}
\Eq(4.13) $$
while, for positive $m\notin\MM_{+}$, an analogous, simpler
expression is found of the form \equ(4.12) with
$$ \AAA_{m,m'} = m^{2} - 4 \| a \|^{2} \d_{m,m'} ,
\Eq(4.14) $$
which are not zero by Lemma 6. Then the following result holds.

\*

\0{\bf Lemma 8.} {\it For $\MM$ chosen according to
Lemma 6, one has for $m\in\MM$
$$ V_{m} = \sum_{m'\in\MM} \DD_{m,m'} [G(v,w)]_{m'} ,
\Eq(4.15) $$
with $\DD$ a $2N \times 2N$ non-singular matrix.}

\*

\0{\it Proof.} It is sufficient to prove that the matrix $\AAA$
with entries \equ(4.13) is not singular. By using \equ(4.3) we can
write the diagonal entries of $\AAA$ as $\AAA_{m,m}=-6a_{m}^{2}$.
Then one realize immediately that one has
$$ \det\AAA = (-1)^{N} \det D_{N}(6,8)
\prod_{m=1}^{N} a_{m}^{2} ,
\Eq(4.16) $$
where $D_{N}(p,q)$ is the $N\times N$ matrix with diagonal
entries $p$ and all off-diagonal entries $q$.
One can easily prove that $\det D_{N}(p,q)=
(p-q)^{N-1}(p+(N-1)q)$. As in our case $p=6$ and $q=8$
(so that $p\neq q$ and $p<(N-1)q$ for all $N\ge 2$)
the assertion follows. \qed

\*

This allows us to extend the analysis of the previous
section to the case in which the function $v_{0}$
is of the form considered here. At the end Theorem 2 is obtained,
with the set $\MM$ chosen according to Lemma 6.

The proof follows the same guidelines sketched in the
previous Sections \secc(2) and \secc(3).
The only differences are that now to each end-point a mode label
$(n_{\vvvv},m_{\vvvv})$, with $m_{\vvvv}\in \MM$ and $n_{\vvvv}=
m_{\vvvv}^{2}$, and an end-point factor $\s_{\vvvv} a_{m_{\vvvv}}$,
with $\s_{\vvvv}={\rm sgn}\,m_{\vvvv}$ are associated. Moreover
the lines $\ell$ with $\g_{\ell}=v$ carry two mode labels
$(n_{\ell},m_{\ell})$ and $(n_{\ell}',m_{\ell}')$, with $m_{\ell}$
and $m_{\ell'}$ both in $\MM$ or in its complement,
and the corresponding propagators are not diagonal any more if
$m_{\ell},m_{\ell}'\in\MM_{+}$, as they are given by $g_{\ell}= \DD_{m,m'}$
(see \equ(4.13) and \equ(4.15)).

Of course the value of $\e_{0}$ depends on the set $\MM$,
and in particular goes to zero when $N\to \infty$
(as $M$ diverges in such a case).

The conclusion is that infinitely many unperturbed solutions
which are trigonometric polynomial with an arbitrary number of
harmonics can be continued in presence of nonlinearities.
The case of polynomials of degree 1 (theorem 1) is the one
usually considered in literature, while the case of polynomials
of higher order (theorem 2) is new. In the latter case
the only requests on the harmonics is that the corresponding
wave numbers have to be close enough to each other
and that larger is their number the larger are their values.

%%%%%%%%%%%%%%%%%%%%%%%%%%%%%%%%%%%%%%%%%%%%%%%%%%%%%%%%%%%%%%%%%%%%%%%%%%
%%%%%%%%%%%%%%%%%%%%%%%%%%%%%%%%%%%%%%%%%%%%%%%%%%%%%%%%%%%%%%%%%%%%%%%%%%
%\*\*
%%%%%%%%%%%%%%%%%%%%%%%%%%%%%%%%%%%%%%%%%%%%%%%%%%%%%%%%%%%%%%%%%%%%%%%%%%
%%%%%%%%%%%%%%%%%%%%%%%%%%%%%%%%%%%%%%%%%%%%%%%%%%%%%%%%%%%%%%%%%%%%%%%%%%
% References
%%%%%%%%%%%%%%%%%%%%%%%%%%%%%%%%%%%%%%%%%%%%%%%%%%%%%%%%%%%%%%%%%%%%%%%%%%
%%%%%%%%%%%%%%%%%%%%%%%%%%%%%%%%%%%%%%%%%%%%%%%%%%%%%%%%%%%%%%%%%%%%%%%%%%

\rife{Bo1}{1}{
J. Bourgain,
{\it Quasi-periodic solutions of Hamiltonian perturbations of 2D
linear Schr\"o\-din\-ger equations},
Ann. of Math. (2)
{\bf 148} (1998), no. 2, 363--439. }
\rife{Bo2}{2}{
J. Bourgain,
{\it Global solutions of nonlinear Schr\"odinger equations},
American Mathematical Society Colloquium Publications, 46.
American Mathematical Society, Providence, RI, 1999. }
\rife{CW1}{3}{
W. Craig and C.E. Wayne,
{\it Newton's method and periodic solutions of nonlinear
wave equations},
Comm. Pure Appl. Math.
{\bf 46} (1993), 1409--1498. }
\rife{GY}{4}{
J. Geng, J. You,
{\it A KAM theorem for one-dimensional Schr\"odinger equation
with periodic boundary conditions},
Preprint, 2004. }
\rife{GM}{5}{
G. Gentile and V. Mastropietro,
{\it Construction of periodic solutions of the nonlinear wave equation
with Dirichlet boundary conditions by the Linsdtedt series method},
J. Math. Pures Appl. (9),
to appear. }
\rife{GMP}{6}{
G. Gentile, V. Mastropietro, M. Procesi,
{\it Periodic solutions for completely resonant nonlinear wave
equations with Dirichlet boundary conditions},
Comm. Math. Phys.,
to appear. }
\rife{KP}{7}{
S.B. Kuksin, J.P\"oschel,
{Invariant Cantor manifolds of quasi-periodic oscillations
for a nonlinear Schr\"odinger equation},
Ann. of Math. (2)
{\bf 143} (1996), no. 1, 149--179. }

\biblio
\ciao